# Moving from Data-Constrained to Data-Enabled Research: Experiences and Challenges in Collecting, Validating and Analyzing Large-Scale e-Commerce Data

Ravi Bapna, Paulo Goes, Ram Gopal and James R. Marsden

*Abstract.* Widespread e-commerce activity on the Internet has led to new opportunities to collect vast amounts of micro-level market and nonmarket data. In this paper we share our experiences in collecting, validating, storing and analyzing large Internet-based data sets in the area of online auctions, music file sharing and online retailer pricing. We demonstrate how such data can advance knowledge by facilitating sharper and more extensive tests of existing theories and by offering observational underpinnings for the development of new theories. Just as experimental economics pushed the frontiers of economic thought by enabling the testing of numerous theories of economic behavior in the environment of a controlled laboratory, we believe that observing, often over extended periods of time, real-world agents participating in market and nonmarket activity on the Internet can lead us to develop and test a variety of new theories. Internet data gathering is not controlled experimentation. We cannot randomly assign participants to treatments or determine event orderings. Internet data gathering does offer potentially large data sets with repeated observation of individual choices and action. In addition, the automated data collection holds promise for greatly reduced cost per observation. Our methods rely on technological advances in automated data collection agents. Significant challenges remain in developing appropriate sampling techniques integrating data from heterogeneous sources in a variety of formats, constructing generalizable processes and understanding legal constraints. Despite these challenges, the early evidence from those who have harvested and analyzed large amounts of e-commerce data points toward a significant leap in our ability to understand the functioning of electronic commerce.

*Key words and phrases:* Large-scale, Internet data, web crawling agents, online auctions, music file sharing.

Ravi Bapna is Associate Professor and Ackerman Scholar, School of Business, University of Connecticut, Storrs, Connecticut 06269, USA e-mail: rbapna@business.uconn.edu. Paulo Goes is Gladstein Professor of Information Technology and Innovation, School of Business, University of Connecticut, Storrs, Connecticut 06269, USA e-mail: paulo.goes@business.uconn.edu. Ram Gopal is GE Capital Endowed Professor of Business, School of Business, University of Connecticut, Storrs, Connecticut 06269, USA e-mail: ram.gopal@business.uconn.edu. James R. Marsden is

# 1. INTRODUCTION AND MOTIVATION

The Internet has spawned an enormous amount of online activity that increasingly permeates our lives and our society. Because of such pervasiveness, the Internet offers data-rich environments that enable new research opportunities. Automated "data collection agents" can gather massive amounts of publicly available data from the Internet, allowing researchers to move from "data-constrained" research to "data-enabled" research. For instance, it is possible to collect information about how online retailers adjust prices over time, what consumers purchase, how consumers search for information, how consumers bid in a given auction or a series of auctions, the degree of interest a news piece generates, what kinds of mp3 files individuals share—all in a very cost-effective and nonintrusive manner. Additional rich sources of web activity data that can necessitate tremendous storage requirements are privately available log files. Generated by the servers that provide web information, these files contain extremely detailed data about each visitor's behavior and actions. But what opportunities and challenges emerge from these enhanced data collection capabilities? Have we reached the stage where we can now actually address micro-level issues using massive micro-level data sets rather than approximate, aggregate-level data? Are electronic markets sufficiently different, or is individual behavior in these markets sufficiently different from previous markets that we need to develop new theories and analyses of the e-markets? What effects do our new data collection capabilities have on our ability to conduct stricter, or perhaps more accurate, tests of existing theories? How does our analysis of this data help us to identify when and how we should be constructing new theories?

Theory formulation often relies heavily on inductive reasoning from repeated observation. Once formulated, theories must be formally tested—tests involving appropriate and adequate data. Not surprisingly, data "limitations" and/or the lack of availability of appropriate data have been common sources of frustration. In economics, the aggregate nature of market data has been a frequently humbling hurdle for those seeking to test theories of consumer behavior and choice. Experimental economists (see, especially, the work of Smith, 1976, 1982, 1991; Plott, 1987; Plott and Sunder 1982, 1988, along with the array of work in Kagel and Roth's handbook, 1995) argued the importance and value of data obtained in controlled laboratory studies using "reward structure(s) to induce prescribed monetary value on actions." In addition to providing aid in theory development and in conducting rigorous empirical tests of economic theory, experimental data provide a data source in situations where relevant field data cannot be obtained (Marsden and Tung, 1999). With all its potential for tracking individual activities, experimental data remain limited because the data are expensive to gather due to the following: (1) the monetary reward mechanisms so critical in Smith's approach, (2) the need to recruit and schedule appropriate subjects, and (3) the experiment development, prototyping and final setup requirements.

In addition to field data and data generated through controlled experimentation, business researchers have often utilized two other data types:

(1) panel data generated by surveying individual consumer or business representatives; and,

(2) simulation data generated by repeated analysis of modeling formulations for varying parameter settings.

Arguably, both of these data types suffer from issues of reliability and appropriateness. Tests may be run on the internal consistency of survey data, but such tests do not address basic reliability linked with subjective nuances and interpretations. Data resulting from simulations may face questions relating to the appropriateness of either the optimization model (does the model accurately represent the intended process?) or the parameter set analyzed (was the range and/or number of cases sufficient?).

The Internet and electronic markets now offer a data source with significant volume and ease of access. But there is much more beyond data volume and access. Users of the Internet operate in an observable venue. In contrast to laboratory experiments, where participants typically work and interact in simulated environments, observing e-commerce activity can provide a direct view into the actions of buyers and sellers, often over extended periods of time. In many situations, whether involved in market or nonmarket activities, individual choices and actions can be tracked and recorded. In this paper, we share our experiences with alternative web-based micro-level e-commerce data from three important domains—online auctions, music file sharing and online retailer pricing.

*Board of Trustees Distinguished Professor, School of Business, University of Connecticut, Storrs, Connecticut 06269, USA e-mail: jmarsden@business.uconn.edu.*

A significant part of this experience has been dealing with statistical challenges including the challenge of structuring appropriate sampling procedures to draw data from millions of occurrences. Consider, for example, the millions (at the time of writing we measured the number of ongoing auctions on eBay to be over 12.5 million; see Section 4 for more details) of online auctions in approximately thirty major categories that run on eBay at any point in time or the millions of files being shared on any of a number of popular peer-to-peer networks. There are temporal issues in capturing data that is often fleeting, for example, a price quote at an online retailer. Additionally, as discussed in Overby (2005), what significance does a $p$-value retain when we are talking about very large sample sizes?

To address such questions and to continue the development of innovative data-capturing software and intelligent automated "data agents," we argue for an interdisciplinary approach bringing together expertise from statistics, computer science and information systems. In doing so, we emphasize the importance of developing and leveraging innovative data gathering and large-scale data analytics as tools to help us improve our understanding of e-commerce activity.

In Section 2 we summarize the attributes of the primary sources of data available to researchers. In Section 3 we dig deeper into three specific ongoing research streams in which we rely on automated Internet-based data collection: online auctions, music file sharing and online retailer pricing. In each case, we summarize the agent-based web crawling and parsing technology developed to collect the data, describe the kinds of new and innovative questions that were addressed using the micro-level observations and present an overview of the key results. In Section 4 we pull together the common methods and technical threads from these research streams and describe the key challenges we see in this new "data-rich" world. Our concluding remarks center on what we view as research arenas that seem likely to move ahead through a bountiful harvest of Internet data.

## 2. DATA TYPES

The value of research data is linked to its ability to aid the researcher in answering interesting research questions. For each research investigation, we seek to identify and utilize the best data type, that is, that data which is most appropriate to help achieve the specific research goals. In what follows, we briefly summarize characteristics of each of four historically common data types and a recently emerging fifth data type (see Hoffman, Marsden and Whinston, 1990).

### 2.1 Data from Historical Observation: Field Data

In summarizing the historical approach of economists, Smith (1987, *The New Palgrave*) argued that economists viewed their field as nonexperimental. Economists historically focused on "field observations"... "seeking to understand the functioning of economies using observations generated by economic outcomes realized over time" (Smith, 1987). Field data provides information on economic outcomes but often suffers from the following undesirable characteristics:

(1) substantial levels of aggregation (e.g., average mortgage rate, average price of gas, total volume sold in geographic region);
(2) researcher removed from data gathering with little or no means to check data accuracy or data gathering procedures; and,
(3) lack of information on how the market outcome was reached (what intermediate actions were taken by economic agents operating in the market?).

Yet field data is often relatively easily and inexpensively obtained and remains the most prevalent data form. Further, field data obtained from accurate company records or regulated markets can be very rich.

### 2.2 Data from Controlled Experiments

For some forty years, Vernon Smith has argued for and repeatedly demonstrated the use of experimental methods to provide a new source of data and a body of "tested behavioral principles that have survived controlled experimental tests" (Smith, 1987). Experimental data is obtained in carefully controlled settings that enable the observation and recording of intermediate actions and choices as well as final market outcomes. In addition, the research is typically the main force in design and completion of the experimentation. Thus, experimental data facilitates avoidance of the three concerns noted above relating to field data. But experimental data is not inexpensive. Costs include experimental design, programming, implementation, identification and scheduling of appropriate subjects, and sufficient monetary rewards for induced value (see Smith, 1976, 1982).

### 2.3 Data from Surveys of Individuals or Business Representatives

Collecting panel data typically involves posing "what if" questions to individuals or business representatives. Respondents are requested to specify their action or decision choice under a variety of scenarios. One difficulty is that the scenarios may not be familiar to respondents. In fact, the scenarios may be far outside any decision-making responsibility or experience of the respondents. Further, the subjects experience no positive reward or negative penalty based upon their response—in Smith's terms, there is no induced value. Unless careful controls are structured and followed, we may have no information on the ability or the incentive of respondents to answer accurately.

Yet, if our research requires information on perceptions, how else can we obtain data except by asking appropriate individuals or business decision-makers?

### 2.4 Data from Simulations

Like controlled laboratory experiments, simulation utilizes an abstraction, a simplification meant to parallel the environment but not recreate that environment. Simulation data is typically generated by repeated model analysis under varying parameter settings. That is, we focus on differences in the outcome space under various alternative specifications of initial conditions and parameter values. Potential difficulties associated with simulation include appropriateness of the model formulation, relevance of the parameter settings utilized and appropriateness of incremental changes in parameter settings. Finally, despite enhanced information technology, simulations can still be time and resource intensive.

Each of the four data types listed above has useful characteristics but also specific limitations. Shugan (2002), in a recent editorial in *Marketing Science*, offered the following rather exceptional summary of data critiques, "Of course, no data source goes unscathed . . . ." But now a fifth data type and source is rapidly evolving. What is it and what can it bring to researchers?

### 2.5 Data from e-Markets and the Internet

As e-markets continue to emerge and develop, they hold promise as a potentially fertile ground for harvesting significant market and nonmarket data. In particular, e-markets can provide a large number of observations on market outcomes and individual behavior in those markets. Internet data gathering is not controlled experimentation. We cannot randomly assign participants to treatments or determine event orderings. Internet data gathering does offer potentially large data sets with repeated observation of individual choices and action. Consider auction markets such as those found on eBay, Sam's Club, Yahoo!, Amazon, and which involve a large number of individual auctions and participants. As illustrated in Bapna, Goes and Gupta (2003a,b), it is quite straightforward to capture all individual bids (including intermediate bids), bid times, bid increments and final outcomes for all bids on a site such as uBid. Contrast this to the limited number of observations and object types (typically common-value object auctions) available in fairly rare field data sets such as that for offshore oil leases in the United States (see Hoffman and Marsden, 1986, and Hoffman, Marsden and Saidi, 1991).

In contrast, participants, both sellers and buyers, in certain electronic activities provide a window on their individual activities—a window from which micro-level data can be recorded, tracked and analyzed. It is the nature of the media that enables this. Researchers such as Smith (1976, 1982, 1991), Plott (1987) and Plott and Sunder (1982, 1988) have included the observing and tracking of microbehavior and individual outcome information in controlled experiments, but the limits on size and scope in this domain clearly separate it from what we are discussing in the Internet domain.

Consider music sharing on peer-to-peer (P2P) networks. An exchange of goods occurs but at no monetary price that is paid to the music "supplier," the individual who permits others to search and access his or her storage media and download files at will. Because of the structure of such networks, it is technically possible to observe and track both individual- and aggregate-level sharing activities. It is also possible to utilize such data in testing theories of digital good (e.g., music) markets, including theories relating to optimal network structure and size (Asvanund et al., 2004), retailer distribution channel selection and consumer behavior.

## 3. INSIGHTS FROM THREE STREAMS OF e-COMMERCE RESEARCH

We offer three specific demonstrations of the collection and utilization of Internet data. In the case of online auctions, we detail how the collection of the entire bidding history of a set of auctions gives us new insights into consumer bidding strategies. For P2P networks, we illustrate how capturing the Internet data enabled an event study analysis of the impact of a specific major legal threat on music sharing activity. Finally, in the context of understanding online retailer pricing strategies, we focus on the process—and the challenges—in capturing and cleaning data from a variety of heterogeneous sources.

Note that this section is not meant to be an exhaustive review of any of the online auction, music piracy or online retailer pricing streams of research. Such an endeavor, while worthy, is outside the scope of this paper. Rather, we have provided three examples from our research streams that illustrate how innovative data collection techniques can enable us to address research questions that had previously eluded rigorous analysis. From each of the three research streams, we focus on one interesting slice intended to demonstrate the breadth of issues and, in doing so, discuss specific characteristics of the data and the process of collection.

### 3.1 Understanding Bidding Strategies in Online Auctions

Auctions have long served as operationally simple resource allocation and pricing mechanisms that coordinate privately held information quickly and accurately to achieve efficient exchange. Auction design has been the focus of significant theoretical (see Rothkopf and Harstad, 1994; McAfee and McMillan, 1987; Milgrom, 1989; Myerson, 1981, for a review) and experimental attention (see Kagel and Roth, 1995, for a review). Traditional (non-Internet) auctions have been the subject of some limited empirical work (Paarsch, 1992; Laffont, Ossard and Vuong, 1995). However, *most research on traditional auctions focuses on the bid-taker's perspective and assumes a certain bidder behavior*. In the traditional, face-to-face auction setting, it might have been reasonable to assume that bidders belonged to a homogeneous, symmetric, risk-neutral group who adopted Bayesian–Nash equilibrium strategies. (Bayesian games are games of incomplete information in which each agent knows her own payoff function but at least one agent exists who is uncertain about another player's payoff function. In the context of a Bayesian game, a Bayesian–Nash equilibrium is one in which each player's course of action is a best response to the other players' strategies.) The advent of eBay and other Internet-based auction sites such as uBid and Sam's Club have brought significant structural and environmental changes to auction markets. These changes led us to investigate the validity of the typical core assumptions in auction theory. The most salient of these, tracing its origin to the fundamental requirement of pursuing a classical game-theoretic analysis, is that the number of bidders is exogenous and is known *ex ante*. However, this assumption is readily violated in online auctions. For example, recent analysis with eBay data reveals that bidder entry is influenced endogenously by the hidden reserve price adopted by the sellers (Bajari and Hortaçsu, 2003). While perhaps tenable in the context of historical face-to-face auctions (where bidders of a certain type may be "invited" to participate in a given auction), the assumptions of a priori known and homogenous bidding strategies as well as a fixed number of bidders quickly break down in most online auctions (Bapna, Goes and Gupta, 2001, 2003a; Bapna, Goes, Gupta and Karuga, 2005). In fact, there is an interesting stream of research in online auctions that focuses on investigating virtually all the traditional issues ranging from revenue equivalence to seller's mechanism design choices (see Bajari and Hortaçsu, 2004, for an exhaustive review).

In addition to revisiting traditional assumptions, online auction research has expanded the auction literature in at least three prominent directions. There is a new emphasis on understanding the heterogeneity of real-world bidders' strategies, a move far away from the vanilla view of auction bidding (Bapna, Goes, Gupta and Jin, 2004; Steinberg and Slavova, 2005). Given the absence of face-to-face interactions, research has begun to center on determining the efficacy of digital reputation systems and on the digitizing of word-of-mouth (see Dellarocas, 2003, as a starting point). There are also burgeoning opportunities in applying new statistical tools, such as functional data and regression analysis, to understand the price formation process of an auction and on utilizing this to create a platform for dynamic mechanism design. Initial work in this direction has looked at formally understanding the general bid-arrival process in online auctions (Shmueli, Russo and Jank, 2004), and at modeling the dynamics of online auctions (Shmueli and Jank, 2006). In large part, these new streams rely heavily on new data extraction and parsing tools that capture complete bidding histories of web-based auctions.

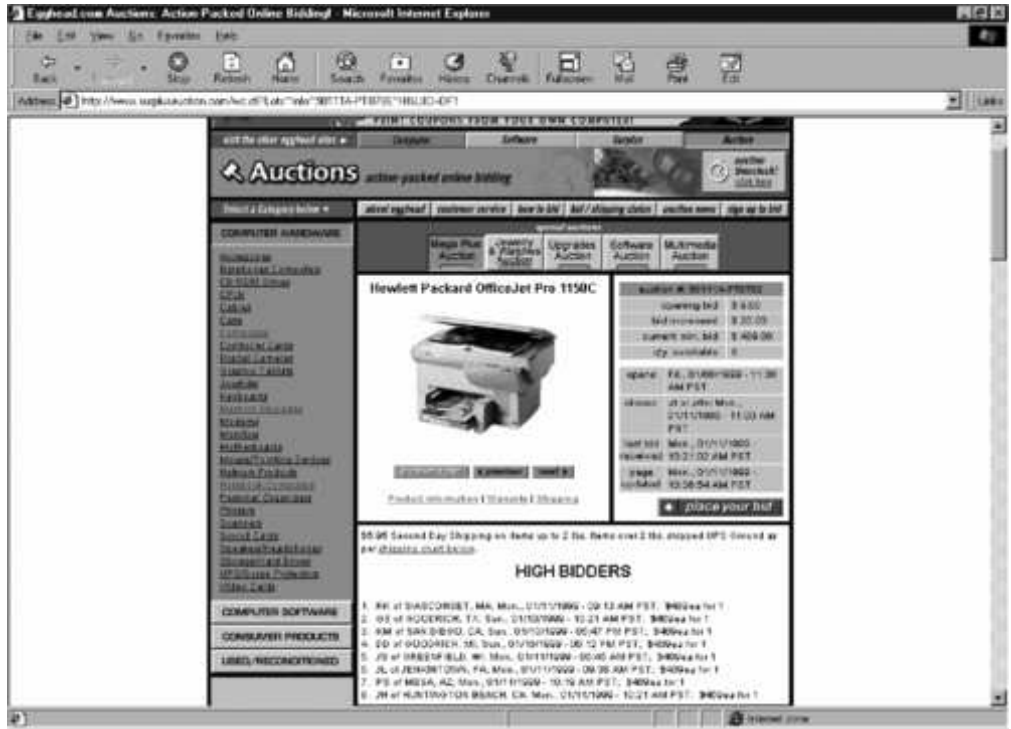

Fig. 1. *An online Yankee auction in progress.*

3.1.1 *Data captured to understand bidding strategies.* Bapna, Goes, Gupta and Jin (2004), hereafter referred to as BGGJ, elaborate on the nature and capture of the data used to analyze heterogeneous bidding strategies in online auctions. They programmed an automated agent to capture information directly from the auction website surplusauction.com, which was later bought by onSale.com. Figure 1 depicts a typical auction in progress, an auction where eight printers are being sold and the two current highest bidders are at $489 while the remaining six are at $469.

The auctions studied by BGGJ were an interesting artifact of the new online environment known as Yankee auctions. In a Yankee auction, there are *multiple identical units* for auction and each auction specifies minimum starting bids and bid increments. At any point, bidders may bid on more than one unit, but all such bids of multiple units must be at the same price. A Yankee auction terminates on or after (most auctions, eBay being the notable exception, have a "going, going, gone" period such that the auction terminates after the closing time has passed and no further bids are received in the last five minutes) a pre-announced closing time. Winning bidders pay the amount they last bid to win the auction. In multi-unit settings, this rule often leads to bidders paying different amounts for the same item. The soft closing time provides a disincentive to last-minute bidding and is designed to attract bids early in the auction. Underlying the web appearance is the actual HTML code that is updated as the auction progresses, as depicted in Figure 2.

It is the HTML documents that are captured in specified small time increments. These documents contain an auction's product description, minimum required bid, lot size and current high bidders. Note that the auction posts only the current winning bidders. Thus to track the progress of the auction and the bidding strategies the same hypertext page has to be captured at frequent intervals until the auction ends. The raw HTML data is then sent to a "parsing" application that searches for specific strings that correspond to the variables of interest, say the price or the lot size. The parsing application condenses the entire information for an auction (including all the submitted bids) into a single spreadsheet. Data is cleansed to screen out auctions in which: (a) sampling loss occurs (due to occasional server breakdowns), and (b) insufficient interest exists (some auctions do not attract any bidders). Sampling loss was detected when there was no intersection between the sets of winning bidders in any two consecutive snapshots of the auction's progress. Table 1 describes the data that was collected and its source.

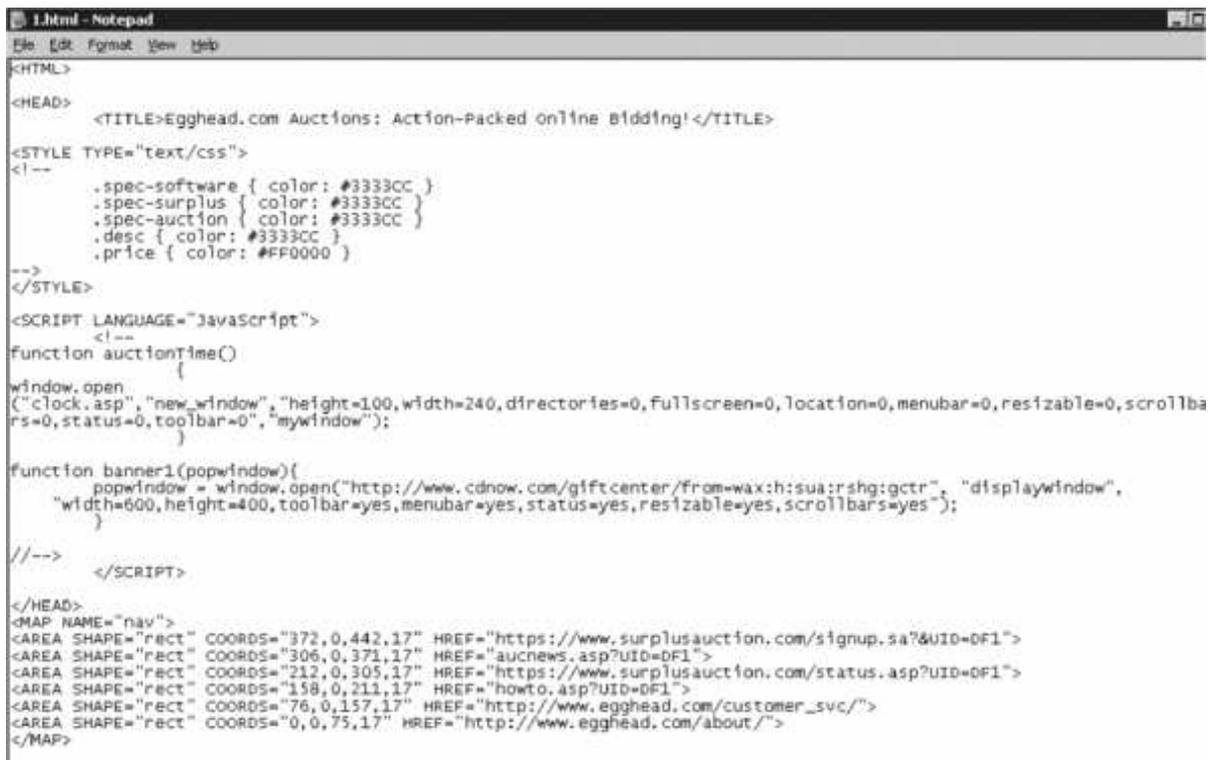

Fig. 2. *Underlying HTML code for a typical Yankee auction.*

BGGJ's data collection spanned a period of two nonconsecutive six-month (plus) periods (in 1999 and 2000) and focused on auctions that sold computer hardware or consumer electronics. There was a one-year gap between the end of the first data collection and the start of the second. Because one of the objectives of the study was to look for any longitudinal shifts in aggregate strategic behavior, the gap in the data collection periods provided the discontinuity required to count experience and exposure as a factor. There were 4580 distinct bidders in this data set. One interesting aspect of the online environment is that, after the one-time registration cost, there is next to no cost for visiting an auction site and placing a ridiculously low bid. Such frivolous bidding could potentially introduce bias in data analysis. BGGJ deemed a bidder frivolous if his final bid was less than 80% of the value of the lowest winning bid. This resulted in 9025 unique bidding data points from 3121 valid bidders participating in 229 auctions. This data was stored in a database to facilitate querying.

3.1.2 *Understanding online bidding strategies.* BGGJ sought to understand how bidding behavior in the online setting differed from the behavior portrayed in the existing literature. They used time of entry, exit and the number of bids as the three key variables in an iterative $k$-means clustering analysis and found a stable taxonomy of bidding behavior as depicted in Figure 3. Bidders can generally be clustered according to those that arrive early and bid multiple times, those that arrive early and bid exactly once, and those that arrive late. The pattern persisted over the two years, with some minor adjustments for the arrival of bidding agent technology.

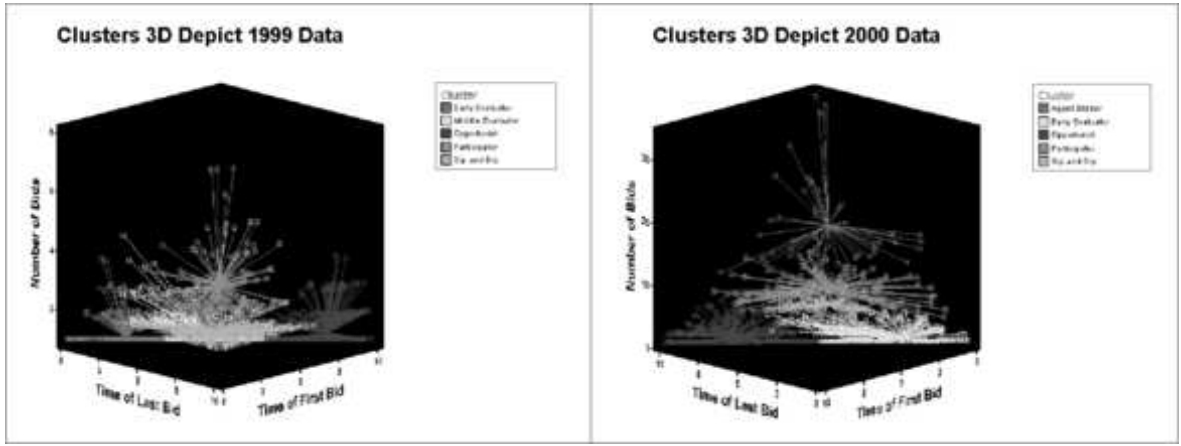

Fig. 3. *Bidder clustering.*

Overall, the classification is robust. In fact, the work of Steinberg and Slavova (2005) was able to replicate the clustering of BGGJ on a different and more recent data set. Steinberg and Slavova (2005) also proposed an alternative classification based on the Yankee auction's allocation rule, one that uses a different variable set (viz. final price, quantity bid and time of entry) than that of BGGJ.

### 3.2 Understanding the impact of music file sharing

While the music recording industry has faced piracy of music ever since the invention of the tape recorder, the level of concern has escalated rapidly with the introduction of new network technologies, including P2P sharing networks. The recording industry and its organization that helps drive legal activities, the Recording Industry Association of America, claim that music sharing activities have had significant negative impacts on the market outcomes for their products. The industry has sought additional regulation and legal penalties, including the threat of individual prosecution, as a means to reduce music sharing. Others (see, e.g., Fader, 2000) have suggested that there may be significant positive effects from P2P sharing activity linked to pre-purchase sampling. What no one had was specific sharing data and sound empirical analysis.

Interestingly, record companies continue to strive to get air-play for their music, even to the extent of continuing *payola*, a "pay-to-play" system of bribery between record companies and radio station personnel. A recent *Business Day* article (August 15, 2005) summarized the continuing payola scandal where New York's Attorney General, Eliot Spitzer, has recently forced SONY BMG to pay $10M to charges of "... bribing radio stations to pump up the airtime allocated to the songs it sells." The article went on:

> Ronald Coase, the Nobel prize-winning economist, explained the practice in 1979. Radio stations own something valuable: songs played more tend to sell more. Competition for airtime develops, but how one conducts the best auction (given that station revenues come primarily from selling audiences to advertisers) is complicated. One view is that radio stations should be faithful to listeners and make choices based only on their disc jockeys' honest musical appreciation. But how do they know what gangsta rap track is top quality?

Table 1
*Operationalization of constructs*

| Construct | Data collected for operationalization | Source |
| --- | --- | --- |
| Product attributes | Type, life cycle, condition | uBid.com and onSale.com |
| Bidder information | Time of entry, exit, number of bids | uBid.com and onSale.com |
| Auction parameters | Lot size, bid increment, duration | uBid.com and onSale.com |
| Bids | Number of bids | uBid.com and onSale.com |
| Market price | Average posted price for a product | Pricegrabber.com |

Payola helps them learn, because record companies will tend to value airtime the most for releases for which they have the highest expectations of sales.

Airtime continues to be precious as a marketing tool. As Spitzer and the commission crack down in the US, competition will perhaps come to entail ever more innovative compensation schemes.

But is not P2P sharing a form of airtime, of increased song play? Could not some significant part of P2P sharing actually provide the prepurchase sampling that record firms seek through payola? Bhattacharjee, Gopal, Lertwachara and Marsden (2006b), hereafter referred to as BGLM, address this very question through the collection and analysis of Internet data. Their data were collected from WinMX, one of the most popular sharing networks that allows users to share, search and download files. BGLM noted that, "On a recent weekend, multiple observations (two each day from Friday to Sunday) indicated an average of 457,475 users sharing 303,731,440 files." Figure 4 presents the results from a search conducted by the authors on the album "Dixie Chicks Home." Figure 4 also includes notation identifying nine variables collected from each search outcome.

BGLM constructed automated tools to capture the rankings released on the weekly Billboard charts and to initiate the searches utilizing this information. It should be noted that once an album appeared on the chart, BGLM continued to include the album in all subsequent searches. Thus the search and the recorded data tend to continually expand. Results from each search were compiled into a relational database.

BGLM's automated search was triggered at random times each day with a random ordering of the albums to be searched. In addition to conducting the searches over all albums after their appearance on the Billboard charts, the authors obtained information on upcoming release dates for 47 albums and began tracking sharing activity on WinMX *prior* to the actual release.

Armed with actual longitudinal sharing data, the authors had the necessary ingredients to analyze and investigate rather than simply "speculate." BGLM summarized their four main empirical findings as follows:

(1) significant piracy opportunity and activity was observed;

(2) the level of sharing opportunities is (*positively*) related to albums' relative chart positions;

(3) there is evidence of both "prepurchase sampling" piracy and "lost-sales" piracy; and,

(4) sharing activity levels provide leading indications of the direction of movement of albums on the Billboard charts.

Significantly, these results hold both for the albums tracked subsequent to their appearance in the top 100 on the Billboard chart and for "new release" albums.

... Though preliminary in nature and covering a limited period, *it is important to emphasize that these are the first results that are based on album-specific actual observations of sharing activity.*

Because the activity of interest occurred on the Internet, it was possible for BGLM to gather the relevant data in real time. As they emphasized, acting solely as observers and not as participants, they were able to track album-specific sharing information across time. Combining this with Billboard chart information and upcoming release dates, they were able to complete a pioneering analysis. The WinMX data represented only part of the Internet data collected by BGLM in their ongoing study of music piracy and the music industry. Their continuing data collection efforts and the related research issues are summarized in Table 2 (see BGLM, 2006a, 2006c for additional details).

### 3.3 Understanding Online Price Dispersion— Heterogeneous Data Sources and Hierarchical Linear Modeling

The widespread adoption of online retailing, coupled with the increased information availability about product, retailer characteristics and nature of competition through the Internet, has renewed research interest in examining how electronic markets operate. Of particular interest is whether the Internet-induced reduction in search costs (Bakos, 1997) does result in the so-called "law of one price" or whether price dispersion persists. The phenomenon of price dispersion occurs when different sellers simultaneously offer to sell the exact same product at substantially different prices. Economic theory suggests that reduced search costs in online markets should lead to reduced price dispersion. At the extreme, all retailers should charge identical prices, equal to the marginal costs of production, for a given product. Even though the Internet was expected to reduce search costs for customers, a series of studies has found online price dispersion to be persistent. Remarkably, the persistence is observed across homogenous goods such as books and CDs (see the variety of products discussed in Bailey, 1998 [Bailey (1998) also included software titles]; Brynjolfsson and Smith, 2000; Lee and Gosain, 2002; Clay and Tay, 2001; Ancarani and Shankar, 2004; Erevelles, Rolland and Srinivasan, 2001). Similar evidence demonstrating existence of price dispersion has also been presented in the service industry (see Clemons, Hann and Hitt, 2002; Bakos et al., 2000).

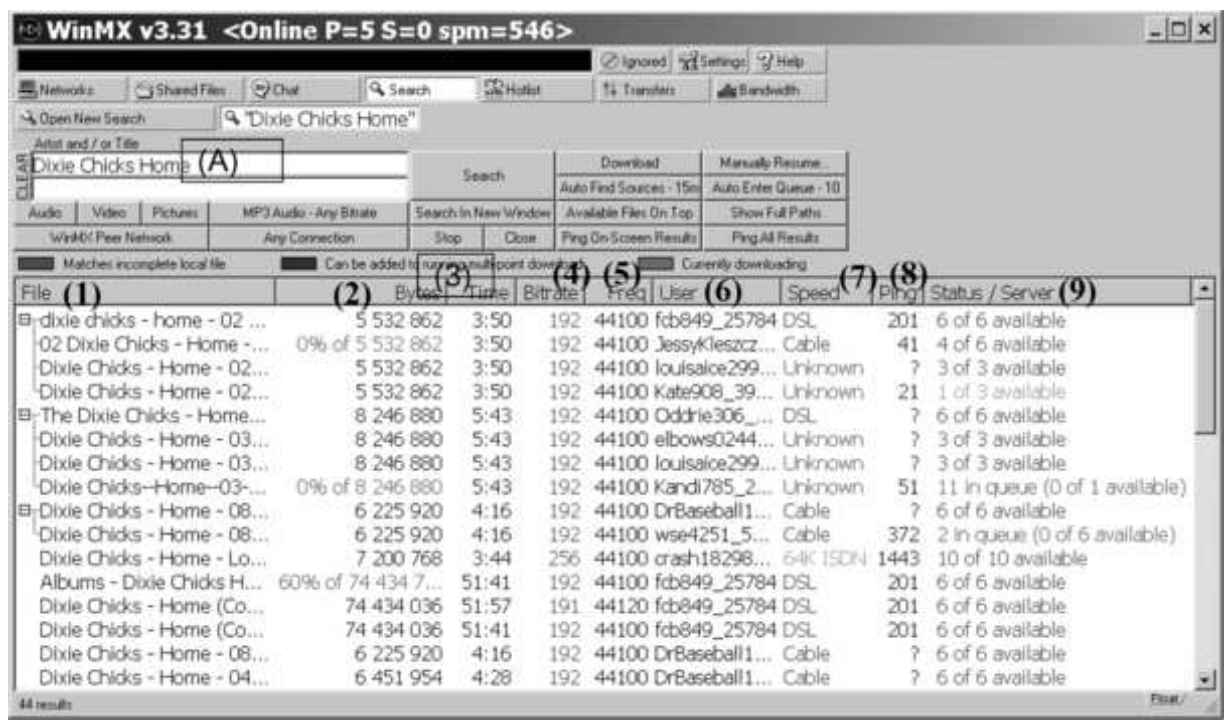

Fig. 4. *A typical WinMX search result.*

Table 2
*Summary of music data collection by BGLM*

| Research problem | Data collected for operationalization | Source |
|---|---|---|
| Does piracy represent prepurchase sampling or lost sales? | Weekly Billboard chart data; daily album-level sharing activity | WinMX |
| What is the impact of piracy on album chart survival? | Weekly Billboard chart data; daily album-level sharing activity (WinMX) | WinMX |
| Are legal threats and/or legal actions effective in deterring online music file sharing? | Daily sharing activity of 2000+ individuals over a two-year period | KaZaA |
| Is prerelease sharing prevalent, and, if so, is it a good predictor of album success? | Daily searches for upcoming release dates, daily album-level sharing activity | www.towerrecords.com www.cdnow.com WinMX |

But why does online price dispersion persist? Potential explanatory variables identified by the literature can be broken into *retailer* and *market* characteristics (see Pan, Ratchford and Shankar, 2003a; Clay, Krishnan and Wolff, 2001). A key challenge in addressing price dispersion issues is that the data needed for analysis exist in a variety of sources and at different levels of aggregation. For instance, retailer characteristics are measured for each retailer, and market characteristics, such as the number of competitors, are measured for each product (see the discussion and approach utilized in Venkatesan, Mehta and Bapna, 2006, and in Raudenbush and Bryk, 2002).

TABLE 3
*Constructs operationalized from a variety of data sources*

| Construct | Data collected for operationalization | Source |
| --- | --- | --- |
| Service quality | Survey ratings obtained by BizRate from online customers. Survey is conducted on a 10-point scale (1 = Poor, 10 = Outstanding) measuring: On Time Delivery, Customer Support, Product Met Expectations and Shop Again | BizRate.com |
| Transactional channels | Dummy coded variables for characterizing the channels through which the retailer offers the products. The channels are first classified as one of the following: pure play (online only), national chain with online presence, and local store(s) with online presence. In addition, we also distinguish retailers that offer mail-order catalog | Manual inspection and verification |
| Size | Rank of retailers based on number of unique visitors to the online store | Alexa.com |
| Competitive intensity | Number of retailers offering an identical product and with service quality ratings available from BizRate | BizRate.com |
| Consumer involvement | Average posted price for a product | BizRate.com |

3.3.1 *Data extraction.* Consider the data collection and cleansing process utilized by Venkatesan, Mehta and Bapna (2006), hereafter referred to as VMB. The online retailing landscape contains a variety of information intermediaries that provide open access to large-sample data on metrics such as comparative product prices, number of competitors, retailer service quality ratings (BizRate.com), and webtraffic (Alexa.com) which can indicate brand or retailer size. Such data are accessible to consumers, to competitors and to researchers. [Interestingly, a recent study, analyzing the representativeness of eight popular shopbots for the online books market, finds significant variance in market coverage of these bots (Allen and Wu, 2003). Their key result is that not all bots are equal and researchers should be judicious in their choice of bots to use for data collection. The eight are Addall, BizRate, Dealtime, ISBN.nu, Mysimon, Pricegrabber, Pricescan and shopping.yahoo.] VMB indicate that while several websites, such as fatwallet.com and resellerratings.com, offer retailer ratings in addition to price quotes, BizRate.com is the only site that follows up with a comprehensive surveying technique in conjunction with retailers. Given its wide coverage in the eight product categories (books, camcorders, DVDs, DVD players, PDAs, printers, scanners and video games) of interest, VMB obtained price and service quality data from BizRate.com. Using BizRate's consumer guide feature, VMB ensured that at least 20 different products were offered in each of the eight categories. In order to have some minimum competitive intensity, they ensured that there were at least seven retailer price quotes for each product. Table 3 provides the summary information about the data collected and each data source. VMB also collected information from Alexa.com. about the size of the retailer measured in terms of number of unique visitors to the retailer's website.

The transaction channel(s) provided by the retailer were coded manually. In order to classify retailer by transaction channels, the retailer's website was inspected to obtain the presence, location and geographical spread of physical store(s). If physical stores were present, and located in more than one state, the retailer was classified as National Brick-n-Click. Retailers with a more localized presence, with physical store locations within one state, were classified as Local Brick-n-Click. If no physical stores were present, the retailer was classified as pure-play. If the retailer offered mail-order catalog services, it was additionally classified as a catalog provider. In rare cases where enough information was not available from the website, VMB searched for news articles regarding the retailer. If that failed to be informative, VMB established telephonic contact with the retailer to obtain the requisite information.

A closer look into the price quotes revealed that the available posted prices were for products in varying condition—new, refurbished or used. Overall, the VMB's web crawling agent collected 22,209 price quotes for 1880 products from 233 retailers. For their data analysis VMB only used prices quoted for items in "New" condition by retailers other than refurb/discounters. VMB excluded retailers with missing service quality ratings (to publish any ratings information on a retailer, BizRate requires a minimum of 30 customer surveys over a period of last ninety days. The published ratings are computed over a rolling window of 90 days), noting that the exclusion of retailers without service ratings was methodologically necessary because service rating was the independent variable in their analysis. Table 4 shows the impact of the bias reduction on the size of the data set suitable for analysis.

Based on this data, analyzed in a hierarchical linear modeling framework, VMB consistently found that:

(1) online retailers who provide better service quality also charge higher prices; and

(2) the influence of service quality on retailer prices varies in a nonlinear fashion with the number of competitors in the market.

Until a certain threshold, retailers with low service quality increase their prices as the number of competitors in the market increases. However, beyond the threshold these retailers decrease their prices with increased competition. In contrast, the high- and medium-service-quality retailers decrease their prices with an increase in the number of competitors until a certain threshold. Beyond this threshold, these high- and medium-service-quality retailers then increase their prices with increased competition. These results provide empirical evidence for theoretical expectations that mixed pricing strategies may exist in online markets (Rosenthal, 1980).

## 4. THE NEED FOR AN INTEGRATED APPROACH

The discussion above illustrates three ongoing research streams that rely on automated real-time collection of Internet data. Albeit in different domains, all three examples share a common sequence of activities with respect to the data. These are (1) automated identification and collection of relevant and related information from different sources; (2) automated cleansing and transformation of data; (3) collating of information at different levels of aggregation; and (4) identification and marking of potentially problematic data for follow-up verification requiring human intervention.

The three examples demonstrate how recent technical advances in collecting Internet-based e-commerce data allow us to pursue a rigorous analysis of seemingly intractable research questions. All three examples involve capturing the online data, integrating the data from multiple sources, cross-validating the information, cleansing the data, offering no threat to the privacy and confidentiality of individuals, and properly storing the data to facilitate appropriate analysis. There are also examples of the use of a similar approach in other research problem domains. For example, Ghose, Smith and Telang (2005) used data collected from Amazon.com to shed new light on the welfare implications of secondary markets for used books. Ghose and Sundararajan (2005), based on software sales data from Amazon.com, empirically quantified the extent of quality degradation associated with software versioning. They also cross-compared their statistical estimation with more subjective web-based data obtained from CNET.com ratings and Amazon.com user reviews. Tang, Montgomery and Smith (2005) employed shopbot-usage data from Media-Metrix and price data from Dealtime.com to demonstrate that an increase in shopbot usage has a positive impact on price dispersion.

All these studies, having extracted large data sets and performed appropriate data analytics, serve as important "proof of concept" examples. Clearly, we are now seeing interesting new results in a variety of problem domains. Yet, we do not have an integrated approach to automated Internet data extraction, validation and storage that enables us to avoid "reinventing the wheel" for every new interesting problem.

We suggest two key steps to move beyond interesting examples and develop a successful Internet data research initiative. First, develop a group of interdisciplinary researchers from three key fields: information systems (IS/MIS), computer science (CS) and statistics (STAT). The data-gathering and storage methods lie mainly within CS and IS/MIS while

ANALYZING LARGE-SCALE E-COMMERCE DATA 13

TABLE 4
*Bias reduction reduces analyzable data*

| Product category | # Posted prices | | # Retailers | | |
|---|---|---|---|---|---|
| | Collected | Analyzed | Collected | Analyzed | # Products |
| Books | 5750 | 2752 | 19 | 9 | 685 |
| Camcorder | 1386 | 882 | 86 | 57 | 57 |
| DVD | 9242 | 5738 | 30 | 15 | 799 |
| DVD player | 547 | 446 | 66 | 51 | 40 |
| PDA | 568 | 479 | 77 | 57 | 32 |
| Printer | 1087 | 906 | 75 | 54 | 36 |
| Scanner | 708 | 574 | 77 | 53 | 31 |
| Video games | 2921 | 1616 | 74 | 42 | 200 |

the STAT community possesses the necessary advanced knowledge in large-scale data analytics. Second, create an environment that can support and foster this interdisciplinary effort. With these two steps in mind, we conceived and developed the Center for Internet Data and Research Intelligence Services (CIDRIS) as a vehicle to initiate an integrated multidisciplinary approach and to work on the development of generalizable techniques. Following the steps suggested above, the initial set of CIDRIS researchers includes faculty from CS, IS/MIS and STAT along with Internet data researchers from a variety of universities (including Minnesota, Arizona, Carnegie Mellon and UT-Dallas). Further, with an initial $400,000 in funding and the top-flight technology of UConn's Gladstein MIS Research Lab and Financial Accelerator, CIDRIS is positioned to nurture and support the necessary interdisciplinary collaborative environment. CIDRIS research is beginning to tackle issues such as the following:

(i) sampling from massive populations of interest (see discussion in Shmueli, Jank and Bapna, 2005). eBay, for example, can have 12,000,000 auctions running simultaneously;

(ii) sampling from dynamic environments (again, see discussion in Shmueli, Jank and Bapna, 2005; Bajari and Hortaçsu, 2004, and Bapna, Jank and Shmueli, 2005)—for example, today's set of products in eBay auctions may be very different from yesterday's;

(iii) integration of automated data extraction tools or wrappers developed by computer scientists (see, e.g., Florescu, Levy and Mendelzon, 1998, and Laender, Ribeiro-Neto, da Silva and Teixeira, 2002)—the major question is, can such online query-focused tools be adapted to facilitate the development of standardized toolkits for data extraction and robust statistical data analytics?

(iv) developing "data validation and analytics toolkits" to address multidisciplinary research questions—can we broaden our approach to include text mining and associated analysis, an arena of potential significant multidisciplinary value? and,

(v) developing data-sharing guidelines to foster academic research—can CIDRIS provide resources to compensate the costs of data gathering while constructing data repositories appropriate for sharing by academic researchers?

But there is a serious caveat to the potential success of CIDRIS and of Internet data collection and analysis in general, a caution linked to legal issues related to automated Internet data gathering. A recent series of legal holdings by a variety of courts provides conflicting outcomes on "trespass to chattels" claims. Winn (2005) offers the following summary and accompanying footnote:

> Once again, while there is considerable uncertainty surrounding the scope of such a (trespass to chattels) claim in light of conflicting case law, the trend in recent cases has been for courts to be more skeptical of such claims and to ask computer owners to tolerate more unwanted interference with the use of computers connected to the Internet. [Trespass was found in Thrifty-Tel, Inc. v. Bezenek, 54 Cal. Rptr. 2d 468 (Cal. Ct. App. 1996); CompuServe, Inc. v. Cyber Promotions, Inc., 962 F. Supp. 1015 (S.D. Ohio 1997); eBay, Inc. v. Bidder's Edge, Inc., 100 F. Supp. 2d 1058 (N.D. Cal. 2000); Register.com, Inc. v. Verio, Inc., 126 F. Supp. 2d 238 (S.D.N.Y.

2000); EF Cultural Travel BV v. Zefer Corp., 318 F.3d 58 (1st Cir. 2003). No trespass was found in Ticketmaster Corp. v. Tickets.com, No. CV99-7654-HLH (VBKx), 2003 U.S. Dist. LEXIS 6483 (C.D. Cal. Mar. 6, 2003); Intel Corp. v. Hamidi, 1 Cal. Rptr. 3d 32 (Cal. Ct. App. 2003); Southwest Airlines Co. v. FareChase, Inc., 318 F. Supp. 2d 435 (N.D. Tex. 2004); Nautical Solutions Mktg. v. Boats.com, No. 8:02-cv-760-T-23TGW, 2004 U.S. Dist. LEXIS 6304 (M.D. Fla. Apr. 1, 2004).]

In a recent *Chicago Bar Association Record* article, Mierzwa (2005) noted the difficulty of showing damages in the case of data "scraping" or automated data gathering and the question of what constitutes acceptance of "terms and conditions" statements. Beyond the issue of damages, there is the potential for a simple denial of service from an Internet site based upon activities such as mining the site too aggressively.

While potential legal issues pose an obstacle to Internet data researchers, they actually provide another important purpose for a center such as CIDRIS to serve as a repository of the latest legal holdings and as an advocate for academic researchers seeking noncommercial, pure research data access. CIDRIS expertise can be directed at developing nondisruptive (or least disruptive) automated data-gathering processes. CIDRIS can work to develop partnering relationships with companies where academic research outcomes are shared with the company while maintaining the right to publish scholarly academic findings without divulging proprietary information. In fact, several of the initial CIDRIS faculty researchers developed and have had such an ongoing relationship with a Fortune 25 company, a partnership now in its fifth year.

As we noted earlier, Internet data gathering is not controlled experimentation. We cannot randomly assign participants to treatments or determine event orderings. Internet data gathering does offer potentially large data sets with repeated observation of individual choices and action.

One additional caveat deserves notice. In many e-commerce domains, identification of participants is limited to "on-line personas," that is, participants self-select id's which keep the participants anonymous to other participants, including "data scrapers." This may limit the usefulness of data in research domains requiring identification of characteristics (demographics, location, etc.) of decision-makers or market participants.

Finally, we note the importance of working with your institution's IT group to establish: (1) that frequent and data-intensive communications will be occurring with specific Internet sites, and (2) the data collection is linked to academic research. This can help avoid a disruptive internal denial of resources.

## 5. CONCLUDING REMARKS

Internet data capture now has the potential to make massive data sets a common reality rather than a rare occurrence. To be successful, we argue the need for an interdisciplinary effort fostered and supported by a center such as CIDRIS. There are technical issues, sampling issues and legal issues. But the availability of large sets of micro-level behavior on the Internet offers the opportunity on a grand scale to pursue the Popperian approach from philosophy of science, that is, structuring stricter and tighter tests of existing theories. It also enables us to address a long-standing problem, aptly articulated by Arthur Conan Doyle through the words he placed on the lips of his most famous character, Sherlock Holmes: "The temptation to form premature theories upon insufficient data is the bane of our profession." Many issues remain, but we are excited about the new world of data-enabled research and the ability to not have to read the following in the data discussions of important research papers:

> "While our analysis is based on a small sample size and the use of numerous proxy variables, it is the best data that we could get. Despite the data issues, the results *appear to be robust.* In future research we hope to address these data issues."

## ACKNOWLEDGMENTS

The research reported here has been supported by the Center for Internet Data Research and Intelligence Services (CIDRIS), Department of Operations and Information Management (OPIM), University of Connecticut. CIDRIS has been funded in part by a 2004 University of Connecticut Provost Grant. More information about CIDRIS is available at cidris.uconn.edu. The authors would like to thank Rob Garfinkel, Kim Marsden and Dave Pingry for their helpful suggestions in improving the manuscript.

# REFERENCES

Allen, G. and Wu, J. (2003). Shopbot market representativeness. In *Proc. International Conference on Electronic Commerce.* Working paper, Dept. Information and Operations Management, Tulane Univ. Available at gallen@tulane.edu.

Ancarani, F. and Shankar, V. (2004). Price levels and price dispersion within and across multiple retailer types: Further evidence and extension. *J. Academy of Marketing Science* **32** 176–187.

Asvanund, A., Clay, K., Krishnan, R. and Smith M. (2004). An empirical analysis of network externalities in peer-to-peer music sharing networks. *Information Systems Research* **15** 155–174.

Bailey, J. (1998). Intermediation and electronic markets: Aggregation and pricing in Internet commerce. Ph.D. dissertation, Technology, Management and Policy, Massachusetts Institute of Technology.

Bajari, P. and Hortaçsu, A. (2003). The winner's curse, reserve prices and endogenous entry: Empirical insights from eBay auctions. *RAND J. Economics* **34** 329–355.

Bajari, P. and Hortaçsu, A. (2004). Economic insights from Internet auctions. *J. Economic Literature* **42** 457–486.

Bakos, J. (1997). Reducing buyer search costs: Implications for electronic marketplaces. *Management Sci.* **43** 1676–1693.

Bakos, J., Lucas, H. C., Jr., Oh, W., Simon, G., Viswanathan, S. and Weber, B. (2000). The impact of electronic commerce on the retail brokerage industry. Working paper, New York Univ.

Bapna, R., Goes, P. and Gupta, A. (2001). Insights and analyses of online auctions. *Comm. ACM* **44**(11) 42–50.

Bapna, R., Goes, P. and Gupta, A. (2003a). Analysis and design of business-to-consumer online auctions. *Management Sci.* **49** 85–101.

Bapna, R., Goes, P. and Gupta, A. (2003b). Replicating online Yankee auctions to analyze auctioneers' and bidders' strategies. *Information Systems Research* **14** 244–268.

Bapna, R., Goes, P., Gupta, A. and Jin, Y. (2004). User heterogeneity and its impact on electronic auction market design: An empirical exploration. *MIS Quarterly* **28** 21–43.

Bapna, R., Goes, P., Gupta, A. and Karuga, G. (2005). Predicting bidders' willingness to pay in online multi-unit ascending auctions. Working paper, Dept. Operations and Information Management, Univ. Connecticut.

Bapna, R., Jank, W. and Shmueli, G. (2005). Consumer surplus in online auctions. Working paper, Dept. Operations and Information Management, Univ. Connecticut. Available at www.sba.uconn.edu/users/rbapna/research.htm.

Bhattacharjee, S., Gopal, R. D., Lertwachara, K. and Marsden, J. R. (2006a). Impact of online technologies on digital goods: Retailer pricing and licensing models in the presence of piracy. *J. Management Information Systems.* To appear.

Bhattacharjee, S., Gopal, R. D., Lertwachara, K. and Marsden, J. R. (2006b). Whatever happened to payola? An empirical analysis of online music sharing. *Decision Support Systems.* To appear.

Bhattacharjee, S., Gopal, R. D., Lertwachara, K. and Marsden, J. R. (2006c). Impact of legal threats on online music sharing activity: An analysis of music industry legal actions. *J. Law and Economics* **49** 91–114.

Brynjolfsson, E. and Smith, M. D. (2000). Frictionless commerce? A comparison of Internet and conventional retailers. *Management Sci.* **46** 563–585.

Clay, K., Krishnan, R. and Wolff, E. (2001). Prices and price dispersion on the web: Evidence from the online book industry. *J. Industrial Economics* **49** 521–539.

Clay, K. and Tay, C. (2001). Cross-country price differentials in the online textbook market. Working paper, Heinz School of Public Policy and Management, Carnegie Mellon Univ.

Clemons, E., Hann, I. and Hitt, L. (2002). Price dispersion and differentiation in online travel: An empirical investigation. *Management Sci.* **48** 534–549.

Dellarocas, C. (2003). The digitization of word-of-mouth: Promise and challenges of online feedback mechanisms. *Management Sci.* **49** 1407–1424.

Erevelles, S., Rolland, E. and Srinivasan, S. (2001). Are prices really lower on the Internet?: An analysis of the vitamin industry. Working paper, Univ. California, Riverside.

Fader, P. S. (2000). Expert report of Peter S. Fader, Ph.D. in record companies and music publishers vs. Napster. United States District Court, Northern District of California.

Florescu, D., Levy, A. and Mendelzon, A. (1998). Database techniques for the World Wide Web: A survey. *SIGMOD Record* **27**(3) 59–74.

Hoffman, E. and Marsden, J. R. (1986). Testing informational assumptions in common value bidding models. *Scandinavian J. Economics* **88** 627–641.

Hoffman, E., Marsden, J. R. and Saidi, R. (1991). Are joint bidding and competitive common value auction markets compatible?—Some evidence from offshore oil auctions. *J. Environmental Economics and Management* **20** 99–112.

Hoffman, E., Marsden, J. R. and Whinston, A. B. (1990). Laboratory experiments and computer simulation: An introduction to the use of experimental and process model data in economic analysis. In *Advances in Behavioral Economics* **2** (L. Green and J. H. Kagel, eds.) 1–30. Ablex, Norwood, NJ.

Ghose, A., Smith, M. and Telang, R. (2005). Internet exchanges for used book. Working paper, New York Univ.

Ghose, A. and Sundararajan, A. (2005). Software versioning and quality degradation? An exploratory study of the evidence. Working paper, New York Univ. Available at papers.ssrn.com/sol3/papers.cfm?abstract_id=786005.

Kagel, J. H. and Roth, A. E., eds. (1995). *The Handbook of Experimental Economics.* Princeton Univ. Press.

Laender, A. H. F., Ribeiro-Neto, B., da Silva, A. S. and Teixeira, J. S. (2002). A brief survey of web data extraction tools. *SIGMOD Record* **31**(2) 84–93.

Laffont, J.-J., Ossard, H. and Vuong, Q. (1995). Econometrics of first-price auctions. *Econometrica* **63** 953–980.

Lee, Z. and Gosain, S. (2002). A longitudinal price comparison for music CDs in electronic and brick-and-mortar markets: Pricing strategies in emergent electronic commerce. *J. Business Strategies* **19** 55–72.

Marsden, J. R. and Lung, Y. A. (1999). The use of information system technology to develop tests on insider trading and asymmetric information. *Management Sci.* **45** 1025–1040.

McAfee, R. P. and McMillan, J. (1987). Auctions and bidding. *J. Economic Literature* **25** 699–738.

Mierzwa, P. (2005). Squeezing new technology into old laws. *CBA Record* **19**.

Milgrom, P. R. (1989). Auctions and bidding: A primer. *J. Economic Perspectives* **3** 3–22.

Myerson, R. B. (1981). Optimal auction design. *Math. Oper. Res.* **6** 58–73. MR0618964

Overby, E. (2005). Size matters: Heteroskedasticity, autocorrelation, and parameter inconstancy in large sample data sets. In *Proc. First Statistical Challenges in E-Commerce Workshop*. Smith School of Business, Univ. Maryland.

Paarsch, H. J. (1992). Deciding between the common and private value paradigms in empirical models of auctions. *J. Econometrics* **51** 191–215.

Pan, X., Ratchford, B. T. and Shankar, V. (2003a). Why aren't the prices of the same item the same at Me.com and You.com?: Drivers of price dispersion among e-tailers. Working paper, Smith School of Business, Univ. Maryland.

Plott, C. R. (1987). Dimensions of parallelism: Some policy applications of experimental methods. In *Laboratory Experimentation in Economics*: *Six Points of View* (A. E. Roth, ed.) 193–219. Cambridge Univ. Press.

Plott, C. R. and Sunder, S. (1982). Efficiency of experimental security markets with insider information: An application of rational-expectations models. *J. Political Economy* **90** 663–698.

Plott, C. R. and Sunder, S. (1988). Rational expectations and the aggregation of diverse information in laboratory security markets. *Econometrica* **56** 1085–1118.

Raudenbush, S.W. and Bryk, A. S. (2002). *Hierarchical Linear Models*: *Applications and Data Analysis Methods*. Sage, Thousand Oaks, CA.

Rosenthal, R. W. (1980). A model in which an increase in the number of sellers leads to a higher price. *Econometrica* **48** 1575–1580.

Rothkopf, M. H. and Harstad, R. M. (1994). Modeling competitive bidding: A critical essay. *Management Sci.* **40** 364–384.

Shmueli, G. and Jank, W. (2006). Modeling the dynamics of online auctions: A modern statistical approach. In *Economics, Information Systems and E-Commerce Research II*: *Advanced Empirical Methods* (R. Kauffman and P. Tallon, eds.). Sharpe, Armonk, NY. To appear.

Shmueli, G., Jank, W. and Bapna, R. (2005). Sampling eCommerce data from the web: Methodological and practical issues. In *ASA Proc. Joint Statistical Meetings* 941–948. Amer. Statist. Assoc., Alexandria, VA.

Shmueli, G., Russo, R. P. and Jank, W. (2004). Modeling bid arrivals in online auctions. Working paper, Smith School of Business, Univ. Maryland. Available at www.smith.umd.edu/ceme/statistics/papers.html.

Shugan, S. M. (2002). In search of data: An editorial. *Marketing Sci.* **21** 369–377.

Smith, V. L. (1976). Experimental economics: Induced value theory. *Amer. Economic Review* **66**(2) 274–279.

Smith, V. L. (1982). Microeconomic systems as an experimental science. *Amer. Economic Review* **72** 923–955.

Smith, V. L. (1987). Experimental methods in economics. In *The New Palgrave*: *A Dictionary of Economics* **2** (J. Eatwell, M. Milgate and P. Newman, eds.) 241–249. Stockton, New York.

Smith, V. L. (1991). *Papers in Experimental Economics*. Cambridge Univ. Press.

Steinberg, R. and Slavova, M. (2005). Empirical investigation of multidimensional types in Yankee auctions. Working paper, Judge Business School, Univ. Cambridge.

Tang, Z., Montgomery, A. and Smith, M. D. (2005). The impact of shopbot use on prices and price dispersion: Evidence from disaggregate data. In *Proc. Workshop on Information Systems and Economics*, Univ. California, Irvine. Available at clients.pixelloom.info/WISE2005/papers/program4.htm.

Venkatesan, R., Mehta, K. and Bapna, R. (2006). Understanding the confluence of retailer characteristics, market characteristics and online pricing strategiess. *Decision Support Systems*. To appear.

Winn, J. K. (2005). Contracting spyware by contract. *Berkeley Technology Law Review* **20** 1345–1359.